\documentclass[12pt]{article}
\usepackage[cmtip,arrow]{xy}
\usepackage{pb-diagram, pb-xy}
\usepackage{amssymb} \usepackage{amsfonts}
\usepackage{amsmath}
\begin{document}

\begin{center}
\textbf{\Large Cohomological rigidity of oriented Hantzsche-Wendt manifolds}\\[0.5cm]
J. Popko, A. Szczepa\'nski\footnote{Corresponding author}\\[0.5cm]

\emph{Institute of Mathematics, University of Gda\'nsk,\\ul.~Wita~Stwosza~57, 80-952 Gda\'nsk, Poland
\footnote{E-mail addresses: matas@univ.gda.pl, jpopko@mat.ug.edu.pl}}
\end{center}\hskip5mm

\date{\today}

\newcommand{\F}{{\mathbb F}}
\newcommand{\Z}{{\mathbb Z}}
\newcommand{\Q}{{\mathbb Q}}
\newcommand{\R}{{\mathbb R}}
\newcommand{\C}{{\mathbb C}}
\newcommand{\h}{{\mathbb H}}
\newcommand{\N}{{\mathbb N}}

\begin{center}
\parbox{12.2cm}{ \footnotesize \textbf{Abstract:}
By Hantzsche-Wendt manifold (for short $HW$-{\em manifold}) we understand any oriented closed Riemannian manifold
of dimension $n$ with a holonomy group $(\Z_2)^{n-1}.$ Two $HW$-{\em manifolds} $M_1$ and $M_2$ are
cohomological rigid if and only if a homeomorphism between $M_1$ and $M_2$ is equivalent to an
isomorphism of graded rings $H^{\ast}(M_1,\F_2)$ and $H^{\ast}(M_2,\F_2).$
We prove that $HW$-{\em manifolds} are cohomological rigid.} 
\end{center}\hskip 5mm

\textbf{MSC2000:} 20H15, 53C29, 57R91, 57S17\\[0.5cm]
\textbf{Keywords:} Hantzsche-Wendt manifold, cohomological rigidity\\[0.5cm]

\section{Introduction}

Let $M^n$ be a flat manifold of dimension $n.$
By definition, this is a compact connected, Riemannian manifold without boundary
with sectional curvature equal to zero. From the theorems of Bieberbach (\cite{Ch}, \cite{S3})
the fundamental group
$\pi_{1}(M^{n}) = \Gamma$ determines a short exact sequence:
\begin{equation}\label{ses}
0 \rightarrow \Z^{n} \rightarrow \Gamma \stackrel{p}\rightarrow
G \rightarrow 0,
\end{equation}
where $\Z^n$ is a torsion free abelian group of rank $n$ and $G$ is a finite group 
which is isomorphic to the holonomy group of
$M^{n}.$ The universal covering of $M^{n}$ is the Euclidean space
$\R^{n}$ and hence $\Gamma$ is isomorphic to a discrete cocompact
subgroup of the isometry group $\operatorname{Isom}(\R^{n}) =
\operatorname{O}(n)\ltimes\R^n = \operatorname{E}(n).$ In the above
short exact sequence $\Z^{n} \cong (\Gamma \cap \R^{n})$ and $p$ can
be considered as the projection $p: \Gamma \rightarrow G\subset
\operatorname{O}(n) \subset \operatorname{E}(n)$ on the first
component.
An orthogonal representation $p$ is equivalent (see \cite{S3}) to a holonomy representation.
That is a homomorphism $\phi_{\Gamma}:G\to \operatorname{GL}(n,\Z),$
given by a formula $\phi_{\Gamma}(g)(z) = \bar{g}z\bar{g}^{-1}$,
where $\bar{g}\in\Gamma, g\in G, z\in\Z^n$ and $p(\bar{g}) = g.$
Conversely, given a short sequence of the form
(\ref{ses}), it is known that the group $\Gamma$ is (isomorphic to)
a Bieberbach group if and only if $\Gamma$ is torsion free.
\vskip 1mm
By Hantzsche-Wendt manifold (for short \emph{$HW$-manifold}) $M^n$ we
understand any oriented flat manifold of dimension $n$ with a
holonomy group $(\Z_2)^{n-1}.$ It is easy to see that $n$ is always
an odd number. Moreover, any HW-manifold has a
diagonal holonomy representation, see \cite{RS}.
It means $\pi_{1}(M^n)$ is generated by $\beta_i =
(B_i,b_i)\in \operatorname{SO}(n)\ltimes\R^n, 1\leq i\leq n,$ where
\begin{equation}\label{genhw}
B_i = \text{diag}(-1,-1,...,-1,\underbrace{1}_i,-1,-1,...,-1)
\end{equation}
and $b_i\in\{0,1/2\}^{n}.$
For other properties of $M^n$ we send a reader to \cite{S3} and to next sections.
We shall need.
\newtheorem{defin}{Definition}
\begin{defin}
{\em (See \cite{KM}.)} Two flat manifolds $M_1$ and $M_2$ are cohomological rigid
if and only if a homeomorphism between $M_1$ and $M_2$ is equivalent to 
an isomorphism of graded rings $H^{\ast}(M_1,\F_2)$ and $H^{\ast}(M_2,\F_2).$
\end{defin}
\vskip 1mm
Our main result is the following theorem. 
\vskip 1mm
\noindent
{\bf Theorem.} {\em Hantzsche-Wendt manifolds are cohomological rigid.}
\vskip 1mm
\noindent
The Theorem answers the question from \cite[problem 4.3]{CMR}. 
\vskip 1mm
\noindent
For the proof we introduce a new presentation of \emph{$HW$-manifolds}.
We consider these manifolds rather as a finite quotient of the torus than a quotient of the $\R^n.$
Here, we use an obvious equivalence $\R^n/\Gamma = (\R^n/\Z^n)/G = T^n/G,$ where $\Gamma$ is a Bieberbach group
from (\ref{ses}). 
According to the definition of $n$-dimensional \emph{$HW$-manifold} we shall define a $(n\times n)$-\emph{$HW$-matrix} $A.$
The analysis of properties of the matrix $A$ is used in the proof.
Moreover, we apply the Lyndon-Hochschild-Serre spectral sequence $\{E_{r}^{p,q}, d_r\}$ of 
the covering $T^n\to T^n/G$
with $\F_2$ coefficients.
Since a holonomy representation $\Phi_{\Gamma}$ is diagonal 
$E_{2}^{p,q} = H^{p}((\Z_2)^{n-1})\otimes H^{q}(\Z^{n}).$
We shall only use the multiplicative structure of the first and second cohomology group.
In particular, we shall consider the properties of the transgression homomorphism
$d_2:H^1(\Z^n)\to H^2((\Z_2)^{n-1}).$  
Finally, another important point of the proof is an isomorphism of
cohomology groups $H^1((\Z_2)^{n-1})$ and $H^1(\Gamma),$ 
which was proved in \cite[Theorem 3.1]{P}.
Hence, we can
consider elements of the image of the transgression homomorphism $d_2$ as
homogeneous polynomials of degree two which are equivalent to polynomial functions. 
\vskip 5mm
\noindent
Let us present a structure of the paper.
In the next section, we give a "new-old" definition of \emph{$HW$-manifold} and we outline the proof of the theorem.
In section three we define \emph{$HW$-matrix} and prove some of its properties.
\vskip 1mm
\noindent 
At the last section, we present 
the proof of the {\bf Main Lemma}.
\section{Proof of the Main Theorem}
Let ${\cal D} = \{g_i\mid i = 0,1,2,3\},$
where $g_i:S^1\to S^1,$ and
$\forall z\in S^1\subset\C,$ 
\begin{equation}\label{dictionary}
g_0(z) = z, g_1(z) = -z, g_2(z) = \bar{z}, g_3(z) = -\bar{z}.
\end{equation}
Equivalently, if $S^1 = \R/\Z, \forall [t]\in \R/\Z,$ 
\begin{equation}\label{dictionary1}
g_0([t]) = [t], g_1([t]) = 
[t+\frac{1}{2}], g_2([t]) = [-t], g_3([t]) = [-t+\frac{1}{2}].
\end{equation}
Let $(t_1,t_2,...,t_n)\in {\cal D}^n$ and $(z_1,z_2,...,z_n)\in T^n = \underbrace{S^1\times S^1\times...\times S^1}_n.$
It is easy to see that ${\cal D} = \Z_2\times\Z_2,$ and $g_3 = g_1 g_2.$ 
For $k=1,2,3$ we have different projections 
\begin{equation}\label{pppp}
p^{(k)}:{\cal D}\to\F_2 = \{0,1\}
\end{equation}
such that 
$p^{(k)}(g_k) = 1$  and for $i = 1,2,..,n$ we have homomorphisms 
\begin{equation}\label{proj}
p^{(k)}\circ pr_i:{\cal D}^n\to {\cal D}\stackrel{p^{(k)}}\to\F_2
\end{equation}
given by the formula
$p^{(k)}\circ pr_i(t_1,t_2,...,t_i,...,t_n) = p^{(k)}(t_i).$
\vskip 2mm
We summing up values of the projections $p^{(2)}$ and $p^{(3)}$ in a table:
\begin{center}
\begin{tabular}{|c|c|c|c|c|}
\hline
& $g_0$ & $g_1$ & $g_2$ & $g_3$\\
\hline
$p^{(2)}$ & $0$ & $1$ & $1$ & $0$\\
\hline
$p^{(3)}$ & $0$ & $1$ & $0$ & $1$\\
\hline
\end{tabular}
\end{center}
\vskip 2mm
\centerline{\tiny{Table 1}}
\vskip 4mm
The next, obvious formula 
\begin{equation}\label{propD}
\forall x\in {\cal D}\hskip 2mm x = p^{(2)}(x)2 +  p^{(3)}(x)3
\end{equation}
will be useful later. 
We can define an action ${\cal D}^n$ on $T^n$ as follows:
\begin{equation}\label{action}
(t_1,t_2,...,t_n)(z_1,z_2,...,z_n) = (t_1 z_1,t_2 z_2,...,t_n z_n).
\end{equation}
We have
\newtheorem{lem}{Lemma}
\newtheorem{prop}{Proposition}
\begin{prop}\label{lemma1}
Let $M^n$ be a HW-manifold of dimension $n.$ Then there exists a subgroup $(\Z_2)^{n-1}\subset {\cal D}^n$
such that $M^n = T^n/(\Z_2)^{n-1},$ where the action $(\Z_2)^{n-1}$ on $T^n$
is defined by {\em (\ref{genhw})} and {\em (\ref{action})}.
\end{prop}
{\bf Proof:}
Let $\pi_1(M^n) = \Gamma$ and
$(B_l,b_l)\in\Gamma$ be the generators (\ref{genhw}), $l = 1,2,..,n.$
On each coordinate, (\ref{dictionary1}) defines $g_j\in {\cal D}, j = 0,1,2,3$ which are determinated
by projections $p^{(1)}\circ pr_{i}, p^{(2)}\circ pr_{i}, p^{(3)}\circ pr_{i}.$
\vskip 2mm
\hskip 120mm $\Box$
\vskip 3mm
\noindent
Let us start to prove that the graded ring $H^{\ast}(M^n,\F_2)$ defines a manifold $M^{n}.$
We have an exact sequence
\begin{equation}\label{hwes}
0\to\Z^{n}\to\Gamma\stackrel{p}\rightarrow(\Z_{2})^{n-1}\to 0,
\end{equation}
where $\Gamma = \pi_1(M^n).$
As we mentioned already in the introduction the image of a holonomy representaion $\Phi_{\Gamma}((\Z_2)^{n-1})$,
is a subgroup of the group of all diagonal matrices of $GL(n,\Z).$ 
Moreover (see \cite{P}) $H^1(\Gamma,\F_2) = (\F_2)^{n-1}$ for any Hantzsche-Wendt group $\Gamma$ of dimension $n.$ 
That is an observation which we shall use during the proof.
\vskip 2mm
\noindent
Since $(\Z_2)^{n-1}\subset {\cal D}^n$ the above 
maps $p^{(k)}\circ pr_{i}, k = 1,2,3$ define homomorphisms from
$(\Z_2)^{n-1}\to\F_2\in Hom((\Z_2)^{n-1},\F_2) = H^1((\Z_2)^{n-1},\F_2) \stackrel{[6]}= H^1(M^n,\F_2).$
Hence we can define elements
$$T_i = (p^{(2)}\circ pr_i)\cup (p^{(3)}\circ pr_i)\in H^2((\Z_2)^{n-1},\F_2),$$
where $\cup$ is a cup product.
It is well known that
$H^{\ast}((\Z_2)^{n-1},\F_2)$ is isomorphic to $\F_{2}[x_1,x_2,...,x_{n-1}].$
Hence the elements
$p^{(k)}\circ pr_i = p^{(k)}_i$ correspond to
\begin{equation}\label{bformula}
\sum_{j=1}^{n-1} p^{(k)}(pr_{i}(b_{j}))x_{j} = \sum_{j=1}^{n-1} p^{(k)}(A_{ji})x_{j}\in \F_{2}[x_{1},x_{2},...,x_{n-1}],
\end{equation}
where $b_1,b_2,...,b_{n-1}$ is the basis of $(\Z_{2})^{n-1}$ and $k = 2,3; i = 1,2,...,n.$
Here the matrix $A_{ij}, i = 1,2,\dots,n-1; j = 1,2,\dots,n$ is related to \emph{$HW$-matrix} 
(Definition \ref{hwm}) from the next section.
\vskip 4mm
\noindent
We shall apply the Lyndon-Hochschild-Serre spectral sequence $\{E_{r}^{p,q}, d_r\}$ of (\ref{hwes}).
Since a holonomy representation $\Phi_{\Gamma}$ is diagonal $E_{2}^{p,q} = H^{p}((\Z_2)^{n-1})\otimes H^{q}(\Z^{n}).$
Hence (see \cite[Corollary 7.2.3 on p. 77]{evens}) we have an
exact sequence (see \cite[p.770]{CMR})
\begin{equation}\label{transg}
H^1(\Z^n,\F_2)\stackrel{d_2}\rightarrow H^2((\Z_2)^{n-1},\F_2)\stackrel{p^{\ast}}\rightarrow H^2(\Gamma,\F_2),
\end{equation}
where $d_2$ is a transgression and $p^{\ast}$ is induced by the above homomorphism $p:\Gamma\to (\Z_2)^{n-1}.$
In what follows we shall prove (see also \cite[Theorem 2.7]{CMR}) that a rank of 
$$\text{Im(}d_2\text{)}\subset H^2((\Z_2)^{n-1},\F_2)\subset H^{\ast}((\Z_2)^{n-1},\F_2)\simeq \F_2[x_1,x_2,...,x_{n-1}]$$
is equal to $n.$ 
\vskip 6mm
\noindent
Let us define a basis $\hat{t_i}, i = 1,2,\dots ,n$ of $H^1(\Z^n,\F_2) =$ Hom$(\Z^n,\F_2).$
For $k\in\Z,$ we shall write $\bar{k} = 0$ if $k$ is even and $\bar{k} = 1$ if $k$ is odd.
Let $(k_1,k_2,\dots,k_n)\in \Z^{n}$ and let 
$$\hat{t}_{i}(k_1,k_2,\dots,k_n) = \bar{k}_i, i=1,2,\dots,n.$$ We have.
\vskip 1mm
\noindent
\begin{prop}\label{rankd2}
{\em $d_2(\hat{t}_i) = T_i = (p^{(2)}\circ pr_i)\cup (p^{(3)}\circ pr_i).$ Moreover elements $T_i, i = 1,2,\dots,n$ are a basis of
{\em $\text{Im(}d_2\text{)}.$}}
\end{prop}
\vskip 1mm
\noindent
{\bf Proof:} 
By Theorem 2.5 (ii) and Proposition 1.3 of \cite{CMR} and using (\ref{bformula}) it follows that
$$d_2(\hat{t}_i) = \sum_{A_{il=1}}x_{i}^{2} + \sum_{i\neq j}x_{i}x_{j},$$
where the second sum is taken for such $i,j$ that 
$$(A_{il},A_{jl})\in\{(1,2),(2,1),(1,3),(3,1),(3,1),(3,2),(2,3)\}.$$
On the other hand 
\vskip 3mm
$$T_l = p^{(2)}_{l}p^{(3)}_{l} = \sum_{i=1}^{n-1}p^{(2)}(A_{il})p^{(3)}(A_{il})x_{i}^{2} +$$ 
\begin{equation}\label{elemmaa}
+ \sum_{1\leq i < j\leq n-1}(p^{(2)}(A_{il})p^{(3)}(A_{jl})
+ p^{(2)}(A_{jl})p^{(3)}(A_{il})x_{i}x_{j}.
\end{equation}
Comparing coefficients of the above two polynomials finishes the proof.  
\vskip 2mm
\hskip 120mm $\Box$
\vskip 3mm
 
The main idea of the proof of rigidity is an application of the above Proposition \ref{rankd2}.
It means, we shaw that any \emph{$HW$-manifold} $M,$ of dimension greater than three, define elements
in the cohomology ring $H^{\ast}(M,\F_2)$ which determines $M$ up to affine equivalence. 
In the {\bf Main Lemma}, we shall prove an existence
of $n$ linear independance elements $T_1, T_2\dots , T_n \in$ Im($d_2$) such that for any $i=1,2,\dots,n \hskip 2mm T_i = p_iq_i.$
At the end of this section we give a method of a reconstruction of \emph{$HW$-group} from the set $\{T_i\}_{i=1,2,\dots,n}.$ 
\vskip 3mm
\noindent
Let us define
\begin{equation}\label{Span}
D = \{y\in\hskip 1mm \text{Im(}d_2\text{)}\mid y\hskip 1mm\text{is a product of two
polynomials of degree}\hskip 1mm 1\}.
\end{equation}
We shall prove that $D$ has less than $n+2$ elements from which we can
reconstruct the basis $T_1, T_2,\dots, T_n$ of Im($d_2$).
\vskip 3mm
\noindent
\vskip 1mm
\noindent
{\bf Main Lemma.}
\noindent
{\em Let $n > 3$, then there are the following possibilities for the structure of the set $D$:
\vskip 1mm
1. $D = \{T_1,T_2,\dots, T_n\};$
\vskip 1mm
2. $D = \{T_1,T_2,\dots, T_n, T_i + T_j\},$ and we can find a polynomial ${\rm p}$ of degree one such that
${\rm p}\mid T_i$ and ${\rm p}\mid T_j$ for some $1\leq i,j\leq n.$
In the second case we can rediscover the set of generators $T_1, T_2,\dots, T_n$.}
\vskip 2mm
\hskip 120mm $\Box$
\vskip 3mm
\noindent
Let $M$ be \emph{$HW$-manifold} of dimension $n.$ 
From the {\bf Main Lemma}, we know that there is a set $D = \{T_1,T_2,...,T_n\}\subset$ Im($d_2$) such that
any $T_i$ is a product of two polynomials $p_i$ and $q_i, i = 1,2,\dots, n$ of a degree one. 
Let $V$ be $(n-1)$-dimensional $\F_2$ vector space. We define a linear map 
$h:V^{\ast}\to {\cal D}^n,$ 
which simple version is (\ref{propD})
such that 
\begin{equation}\label{rediscover}
h_i(x) = p_i(x)2 + q_i(x)3,\hskip 1mm \text{for}\hskip 2mm i = 1,2,\dots,n, 
\end{equation}
where $p_i, q_i \in V\simeq V^{\ast\ast}.$
Hence, through formulas 
(\ref{bformula}), (\ref{elemmaa}) and the Table 1, Im$(h),$ defines a Hantzsche-Wendt group.
\vskip 1mm
\noindent
\newtheorem{examp}{Example}
\begin{examp} 
1. Let $V = \text{gen}\{x_1,x_2,x_3\}$ and $D = \{x_1^2 + x_{1}x_{2}, x_{1}x_{2} + x_{1}x_{3} + x_2^2 + x_{2}x_{3}\}.$
Put $p_1 = x_1, q_1 = x_1 + x_2, p_2 = x_1 + x_2, q_2 = x_2 + x_3.$
Hence a homomorphism $h(x_{1}^{\ast}) = (1,2), h(x_{2}^{\ast}) = (3,1)$ and $h(x_{3}^{\ast}) = (0,3).$
Here $x_{1}^{\ast}, x_{2}^{\ast}, x_{3}^{\ast}$ is a dual basis of $V^{\ast}.$
Finally we define a subgroup of ${\cal D}^2$ which generators are rows of the matrix
$$
\left[
\begin{array}[c]{cc}
1 & 2\\
3 & 1\\
0 & 3 
\end{array}
\right].
$$
\vskip 2mm
\noindent
2. Let $\Z_{2}^{n-1}\subset {\cal D}^n$ be a \emph{$HW$-group}, and $D$ a set from the Proposition \ref{rankd2}.
Assume that $D = \{p_{1}q_{1}, p_{2}q_{2},\dots,p_{n}q_{n}\}.$ Then
$$h_i(x) = p_{i}(x)2 + q_{i}(x)3 = p(x_i)2 + q(x_i)3 = x_i.$$
Hence for $x\in\Z_{2}^{n-1}, h(x) = x$ and {\em Im(}$h${\em )} = $\Z_{2}^{n-1}.$
\end{examp}
Let $\phi:H^{\ast}(M_1,\F_2)\to H^{\ast}(M_2,\F_2)$ be an isomorphism of cohomology rings of \emph{$HW$-manifolds}
$M_1$ and $M_2.$ From the {\bf Main Lemma} for the both manifolds we have the sets od elements $D_1$ and $D_2$
such that $\phi(D_1) = D_2.$ Hence we obtain the affine equivalence manifolds $M_1$ and $M_2.$ 
\vskip 2mm
\hskip 120mm $\Box$
\vskip 3mm
\section{Properties of Hantzsche-Wendt matrices}
Let us illustrate the Proposition \ref{lemma1} on two \emph{$HW$-manifolds} of dimension 5 , (see \cite{S3}). 
We shall denote by $\Gamma_1$ and $\Gamma_2$ its fundamental groups. 
\begin{examp}\label{example1}
A group $\Gamma_1\subset E(5)$ is generated by 
$$(B_1,(1/2,1/2,0,0,0)), (B_2,(0,1/2,1/2,0,0)),$$ 
$$(B_3, (0,0,1/2,1/2,0)), (B_4, (0,0,0,1/2,1/2)).$$
From above $\R^5/\Gamma_1\simeq T^5/(\Z_2)^4,$ where
$(\Z_2)^4\subset {\cal D}^5$ is defined by  
$$(g_1,g_3,g_2,g_2,g_2), (g_2,g_1,g_3,g_2,g_2),$$
$$(g_2,g_2,g_1,g_3,g_2), (g_2,g_2,g_2,g_1,g_3).$$
\vskip 1mm
\noindent
Moreover a group $\Gamma_2\subset E(5)$ is generated by 
$$(B_1,(1/2,0,1/2,1/2,0)), (B_2, (0,1/2,1/2,1/2,1/2)),$$
$$(B_3,(1/2,1/2,1/2,1/2,1/2)), (B_4,(1/2,0,1/2,1/2,1/2)).$$
Hence, $\R^5/\Gamma_2\simeq T^5/(\Z_2)^4$ where 
generators of a group $(\Z_2)^4\subset {\cal D}^5$ are following 
$$(g_1,g_2,g_3,g_3,g_2), (g_2,g_1,g_3,g_3,g_3),$$ 
$$(g_3,g_3,g_1,g_3,g_3), (g_3,g_2,g_3,g_1,g_3).$$
\end{examp}
In what follows we shall write $i$ for $g_i$, $i = 0,1,2,3.$
Let $A$ be a $(n\times m)$ matrix with coefficients $A_{ij}\in {\cal D}.$ For short $A\in {\cal D}^{n\times m}.$ 
Let $A_i$ ($A^j$) denote $i$-row ($j$-column) of a matrix $A.$
\begin{defin}\label{hwm}
By \emph{$HW$-matrix} we shall understand a matrix $A\in {\cal D}^{n\times n}$ such that $A_{ii} = 1$, $A_{ij}\in\{2,3\}$ for $i\neq j, 1\leq i,j\leq n$
and if $X\subset\{1,2,...,n\}$ and $1\leq \# X\leq n-1$ then the row $\sum_{i\in X} A_i$ has $1$ on a some position.
\end{defin} 
\begin{lem}\label{lemma2}
Any HW-manifold of dimension $n$ defines a $(n\times n)$ \emph{$HW$-matrix}. 
\end{lem}  
{\bf Proof:} 
Let $(\beta_i, b_i), 1\leq i\leq n-1$ be generators of the fundamental group of some $n$-dimensional HW-manifold $M.$
Then $i$ - generator defines $i$-row of some $(n\times n)$ \emph{$HW$-matrix}, cf. (\ref{genhw}), (\ref{dictionary1}).
See also Example \ref{example1} and Proposition \ref{lemma1}. The last row is defined by the product $\beta_{1}\beta_{2}\dots\beta_{n-1}$
or equivalently is a sum of the first $(n-1)$ rows.  
It is easy to see that the first property of the above matrix follows from a definition, see \cite[p. 4]{MR}.  
Since a holonomy group $(\Z_2)^{n-1}$ acts free on $T^n$ (or equivalently $\pi_{1}(M)$ is a torsion free group) the last part of lemma follows.  
\vskip 2mm
\hskip 120mm $\Box$
\vskip 3mm
We shall present some properties of \emph{$HW$-matrices}.
\newtheorem{rem}{Remark}
\begin{rem}\label{remark1}
Let $\sigma\in S_n$ and let $P_{\sigma}$ be the corresponding permutation matrix.
It is not difficult to see that if $A$ is \emph{$HW$-matrix} then $P_{\sigma}AP_{\sigma}^{-1}$
also satisfies conditions of the Definition \ref{hwm}.
Moreover, if $A'$ is a conjugation matrix of $A,$ where conjugation means exchange at
some column numbers $2$ for $3,$ then $A'$ is again a \emph{$HW$-matrix}. 
The \emph{$HW$-matrix} is related to the matrix defined on page 6 of {\em \cite{MR}}.
\end{rem}
\begin{rem}
Let $A$ be a $(n\times n)$ \emph{$HW$-matrix}. 
Then  
\vskip 2mm
\begin{equation}\label{matrix1}
(p^{(2)}+p^{(3)})(A) = \left[
\begin{smallmatrix}
0&1&1&\ldots &1&1\\
1&0&1&\ldots &1&1\\
\vdots & \vdots &\vdots &\ddots &\vdots &\vdots\\
1&1&\dots &1&0&1\\
1&1&\ldots&1&1&0
\end{smallmatrix}\right].
\end{equation}
\vskip 2mm 
\end{rem}
\vskip 1mm
\noindent
Let $A\in {\cal D}^{m\times n}$ be a $(m\times n)$ matrix with coefficients in ${\cal D}$ and
$(\alpha_1,\alpha_2,...,\alpha_n)\in \{2,3\}^n.$ By $p^{(\alpha)}(A)$ we shall
understand a $(m\times n)$-matrix with coefficients in $\F_2$ which a $i$-column is equal to
$p^{(\alpha_i)}(A^i).$ 
\vskip 1mm
\noindent
Let $M$ be a matrix. By defect of $M$ we shall understand a number 
$$d(M) = \{\text{number of columns of M}\} - \text{rk}(M).$$
\begin{lem}\label{defect}
1. Let $M_1$ be a matrix $M$ from which we remove some columns. Then 
$$d(M_1)\leq d(M),$$
\vskip 1mm
\noindent
2. If $A$ is a \emph{$HW$-matrix} of dimension $n$ and $\alpha\in\{2,3\}$, then
$$d(p^{(\alpha)}(A))\leq 1.$$
\end{lem}
{\bf Proof:} The first statment is clear. For the proof of a second one, let us assume
that $d(p^{(\alpha)}(A)) > 1.$ 
Hence $rk(p^{(\alpha)}(A)) < n-1.$ By definition there exists a non-trivial $X\subset\{1,2,...,n-1\},$
such that $\sum_{i\in X}p^{(\alpha)}(A_i) = 0.$
Finally $p^{(\alpha)}(\sum_{i\in X}A_i) = \sum_{i\in X}p^{(\alpha)}(A_i) = 0.$
This contradicts the definition \ref{hwm}.
\vskip 2mm
\hskip 120mm $\Box$
\vskip 3mm
\begin{lem}\label{defect1}
Let $m < n$ and $W\in {\cal D}^{m\times n}$ is a sub-matrix of some $(n\times n)$ \emph{$HW$-matrix}.
Then $rk(p^{(\alpha)}(W)) = m.$
\end{lem}
{\bf Proof:} Similar to the proof of the last Lemma.
\vskip 2mm
\hskip 120mm $\Box$
\vskip 3mm
A symmetric $(m\times m)$ matrix $A\in (\F_2)^{m\times m}$ defines a nonoriented graph, graph$(A)$ with
set of vertices $\{1,2,...,m\}$ and two different vertices $i$ and $j$ are
connected if and only if $A_{ij} = 1.$ We say that a matrix $A$ is connected if a graph$(A)$ is connected.
Let $A\in {\cal D}^{m\times m}$ be a symmetric matrix, then $p^{(i)}(A)$ are symmetric with
coefficientes in $\F_2, i=2,3.$
We shall write $i\sim_{2} j$ if $i,j$ are at the same connected component of a matrix $p^{(2)}(A).$
Similar definition is for a relation $i\sim_{3} j.$
\vskip 1mm
\begin{lem}
Let a \emph{$HW$-matrix} $M$ have the following decomposition on the blocks:
\begin{equation}\label{matrix2}
M = \left[
\begin{smallmatrix}
\ast & 2 & \ast\\
C & A & D\\
\ast & 3 & \ast
\end{smallmatrix}\right],
\end{equation}
where $A$ is a symmetric matrix and $2,3$ are block matrices with all rows and columns equal $2$ and $3$ correspondingly. Then
\vskip 1mm
\noindent
{\em (I)} if $i\sim_2 j\hskip 2mm\Longrightarrow\hskip 2mm D_i = D_j;$
\vskip 1mm
\noindent
{\em (II)} if $i\sim_3 j\hskip 2mm\Longrightarrow\hskip 2mm C_i = C_j.$ 
\end{lem}
{\bf Proof:} For the proof of (I) let us assume that $i,j$ (where $i < j$) are connected by a 2-edge; i.e $A_{i,j} = 2.$
Let $r$ be some column of a matrix $D.$ Let us consider a diagonal submatrix of the matrix $M$ related to $(i,j,r)$.
It looks like
\begin{equation}\label{matrix3}
\left[
\begin{smallmatrix}
1 & 2 & a\\
2 & 1 & b\\
3 & 3 & 1
\end{smallmatrix}\right].
\end{equation}
The sums of the first two columns are zero. Since a Lemma \ref{defect1} a sum of elements of the last one
is not zero. Hence $a = b.$ We have just proved that if $A_{i,j} = 2$ then $D_i = D_{j}.$
It also means that if $i\sim_2 j$ then $D_i = D_j.$
The proof of the second point of the lemma is similar. 
\vskip 2mm
\hskip 120mm $\Box$
\vskip 3mm 
\noindent
The next lemmas are about possibilites of complement of some matrices to a \emph{$HW$-matrix}.
We shall first consider an odd case.
\begin{lem}\label{oddlemma}
Let $A\in {\cal D}^{m\times m}$ be a symmetric matrix with $1$ on the diagonal and $\{2,3\}$
off the diagonal with a column sums equal to 1. Assume that $m > 1.$ Then a matrix
\begin{equation}\label{matrix4}
K_A = \left[
\begin{smallmatrix}
2\\
A\\
3
\end{smallmatrix}\right],
\end{equation}
cannot be complement to \emph{$HW$-matrix}.
\end{lem}
{\bf Proof:} 
Let us assume that there axists a \emph{$HW$-matrix}
\begin{equation}
\left[
\begin{smallmatrix}
\ast & 2 & \ast\\
C & A & D\\
\ast & 3 & \ast
\end{smallmatrix}\right].
\end{equation}
From assumption $m$ is an odd number and hights of the blocks $2$ and $3$ are also odd.
We shall use induction. For $m = 3$
\begin{equation}
A = \left[
\begin{smallmatrix}
1 & a & a\\
a & 1 & a\\
a & a & 1
\end{smallmatrix}\right].
\end{equation}
Here $a = 2$ or $3.$
If $a = 3$ then rk$(p^{(2)}(A)) = 1$ and $d(p^{(2)}(A)) = 3-1 = 2 > 1.$
From Lemma \ref{defect} it is impossible. For $a = 3$ the proof is the same.
Let us assume that $m > 3.$
\vskip 2mm
\noindent
1. We shall consider a matrix $p^{(2)}(A).$ We claim that there is no such decomposition as
$$p^{(2)}(A) = B\oplus E,$$ such that a dimension of a matrix $B$ is odd and $> 1.$
In fact, in that case
\begin{equation}
A = \left[
\begin{smallmatrix}
\tilde{B} & 3 \\
3 & \tilde{E}
\end{smallmatrix}\right].
\end{equation}
Since a column sums of $A$ are equal to $1$ and height of a block $3$ under $\tilde{B}$ is even, a column sums of $\tilde{B}$ are $1.$
If $K_A$ has complement then $K_{\tilde{B}}$ has a complement (where a dimension of a block $3$ is greater on a dimension of $E$). But 
by induction it is impossible, since $1 <$ dimension$(\tilde{B}) < m.$
\vskip 2mm
\noindent
2. We claim that there is no such a nontrivial decomposition as
$$p^{(2)}(A) = B\oplus E\oplus F.$$
\vskip 1mm
\noindent
In fact since $m$ is odd we have two possibilities:
\vskip 1mm
(a) dimension of one component is odd and other components have dimension even
\vskip 1mm
(b) dimension of all components are odd.
\vskip 1mm
\noindent
In the case (a) dim$(B\oplus E) > 1$ and odd. Hence we consider decomposition
$p^{(2)}(A) = (B\oplus E)\oplus F.$ But it is a previous case 1.
\vskip 1mm
\noindent
In case (b), since $m > 3$ there exists a component (for example $B$) which dimension is $> 1.$
In that case we have a decomposition
$p^{(2)}(A) = B\oplus (E\oplus F)$ which was already considered in the point 1.
\vskip 2mm
\noindent
3. By definition we have a decomposition
$$p^{(2)}(A) = B_1\oplus...\oplus B_s,$$
where all components are connected matrices.
From the above we can assume that $s = 2$ and odd component has a graph equal to one point or $s=1.$
Equivalently, 
\vskip 2mm
\noindent
(a)\hskip 5mm
$A = \left[
\begin{smallmatrix}
1 & 3 \\
3 & B
\end{smallmatrix}\right]$
and $p^{(2)}(B)$ is connected or 
\vskip 2mm
\noindent
(b)\hskip 5mm $p^{(2)}(A)$ is connected.
\vskip 2mm
\noindent
In the first case
\begin{equation}
p^{(3)}(A) = \left[
\begin{smallmatrix}
1 & 1 \\
1 & p^{(3)}(B) 
\end{smallmatrix}\right].
\end{equation}
Hence $p^{(3)}(A)$ is connected. Summing up, we have 
\vskip 2mm
\noindent
(a)\hskip 5mm
$A = \left[
\begin{smallmatrix}
1 & 3 \\
3 & B
\end{smallmatrix}\right]$
and both $p^{(2)}(B)$ and $p^{(3)}(A)$ are connected or
\vskip 2mm
\noindent
(b)\hskip 5mm
$p^{(2)}(A)$ is connected. 
\vskip 5mm
\noindent
If we exchange $p^{(2)}$ for $p^{(3)}$ in the above points 1., 2. and 3. with the similar arguments, we obtain finally two cases:
\vskip 2mm
\noindent
(a)\hskip 5mm
$A = \left[
\begin{smallmatrix}
1 & 3 \\
3 & B
\end{smallmatrix}\right]$
and both $p^{(2)}(B)$ and $p^{(3)}(A)$ are connected or
\vskip 1mm
\noindent
(b)\hskip 5mm both $p^{(2)}(A)$ and $p^{(3)}(A)$ are connected.
\vskip 2mm
We come back to the beginning of the proof. We shall try to figure out matrices $C$ and $D.$
From definition of $\sim_{3}$ and because $p^{(3)}(A)$ is connected we conclude that all rows of the matrix $C$
are identical. By conjugation 
we can assume that $C = 2.$
Using the same arguments and definiton of  $\sim_{2}$ together with a connectedness of $p^{(2)}(B)$ we
conclude that with exception of the first row, all rows of the matrix $D$ are the same.
By conjugation and permutation we can assume that the first row of the matrix $D$
is equal to $[2,...2,3,...,3].$ 
All other rows of a matrix $D$ consist only $3.$
Summing up a matrix 
$$W = [C\hskip 2mm A\hskip 2mm D]$$ is following
\begin{equation}
\left[\begin{array}{c}2\\2\end{array}
\begin{bmatrix}1&3\\ 3 & B\end{bmatrix}
\begin{array}{c}2\\3\end{array}
\begin{array}{c}3\\3\end{array}\right].
\end{equation}
Apply homomorphisms: $p^{(3)}, [p^{(2)}, p^{(3)}], p^{(2)}, p^{(2)}$ to the corresponding columns we get a matrix
\begin{equation}
W' = \left[\begin{array}{c}0\\0\end{array}
\begin{bmatrix}1&1\\0& p^{(3)}(B)\end{bmatrix}
\begin{array}{c}1\\0\end{array}
\begin{array}{c}0\\0\end{array}\right].
\end{equation}
We have rk$W' = 1 +$ rk$(p^{(3)}(B)).$
From assumption a sums of columns of a matrix $A$ are equal to $1$. Hence
a sums of columns of a matrix 
\begin{equation}
(p^{(2)}, p^{(3)})A = \left[
\begin{smallmatrix}
1 & 1 \\
0 & p^{(3)}(B)
\end{smallmatrix}\right]
\end{equation}
are also equal to $1$ and a sums of columns of a matrix $p^{(3)}(B)$ are equal to $0.$ It means rk($p^{(3)}(B)) < m - 1$ and also
rk$(W') < m.$ 
From Lemma \ref{defect1} 
$$\text{rk}(W') = \text{rk}(W) = \text{number of rows}\hskip 1mm(W) = m.$$
Hence a matrix $W$ cannot be a matrix of some rows of \emph{$HW$-matrix}.
\vskip 2mm
\noindent
We have to still consider a case when matrices $p^{(2)}(A)$ and  $p^{(3)}(A)$ are connected.
Similar to the above consideration, using relation $\sim_2$ and $\sim_3$ plus conjugation we can assume
that 
$$[C\hskip 2mm A\hskip 2mm D] = [2\hskip 2mm A\hskip 2mm 3].$$
Hence all nonempty sums of rows of a matrix $A$ include $1.$ For $m > 1$ it is impossible.
\vskip 2mm
\hskip 120mm $\Box$
\vskip 3mm
\noindent
The next lemma is an even version of the Lemma \ref{oddlemma}.
\begin{lem}\label{evenlemma}
Let $A\in {\cal D}^{m\times m}$ be a symmetric matrix with $1$ on the diagonal and $\{2,3\}$
off the diagonal with a column sums equal to $3.$ Assume that $m > 1.$ Then a matrix
\begin{equation}
K_A = \left[
\begin{smallmatrix}
2\\
A\\
3
\end{smallmatrix}\right],
\end{equation}
cannot be a complement to some \emph{$HW$-matrix}.
\end{lem}
{\bf Proof:} As in the proof of the previous lemma let us assume that there exists a \emph{$HW$-matrix}
\begin{equation}
\left[
\begin{smallmatrix}
\ast & 2 & \ast\\
C & A & D\\
\ast & 3 & \ast
\end{smallmatrix}\right].
\end{equation}
From assumption and Definition \ref{hwm} $m$ is an even number and a hight of the block $2$ is even and $3$ is odd.
We shall use induction. For $m = 4.$
\vskip 1mm
\noindent
1. On the beginning let us consider the case, where $p^{(2)}(A)$ is not connected.
We have two cases of matrices of dimension $4$:
\vskip 2mm
\noindent
(a)\hskip 5mm
$A = \left[
\begin{smallmatrix}
1 & 3\\
3 & B
\end{smallmatrix}\right],$
where $B$ has a dimension $3$ and
\vskip 2mm
\noindent
(b)\hskip 5mm
$A = \left[
\begin{smallmatrix}
B & 3\\
3 & E
\end{smallmatrix}\right],$ where matrices $A, B$ have rank two. 
\vskip 2mm
\noindent
The case (a) is impossible since $1+3\neq 3.$
In the case (b) matrices $A$ and $B$ are symmetric with columns sums equal to $3.$
Hence $B = E = \left[
\begin{smallmatrix}
1 & 2\\
2 & 1
\end{smallmatrix}\right],$ and
\begin{equation}
p^{(2)}(A) = \left[
\begin{smallmatrix}
1 & 1 & 0 & 0\\
1 & 1 & 0 & 0\\
0 & 0 & 1 & 1\\
0 & 0 & 1 & 1
\end{smallmatrix}\right].
\end{equation}
\vskip 1mm
\noindent
From the other side a matrix $p^{(2)}(K_{A})$ has rows of $1$ ($p^{(2)}(2) = 1$) and rows of $0$ ($p^{(2)}(3) = 0$).
These rows are linear combination of rows of $p^{(2)}(A)$ and
$$\text{rk}p^{(2)}(K_A) = \text{rk}p^{(2)}(A) = 2.$$
Finally d($K_A$) = $4-2 = 2 > 1$ and from Lemma \ref{defect} we are done.
\vskip 1mm
\noindent
2. As the second step let us consider the case where $p^{(3)}(A)$ is not connected.
We have to consider two cases of matrices of dimension $4$:
\vskip 2mm
\noindent
(a)\hskip 5mm 
$A = \left[
\begin{smallmatrix}
1 & 2\\
2 & B
\end{smallmatrix}\right],$ and
\vskip 2mm
\noindent
(b)\hskip 5mm
$A = \left[
\begin{smallmatrix}
B & 2\\
2 & E
\end{smallmatrix}\right],$ and $B$ and $E$ have dimension $2.$
\vskip 2mm
In the case (a) a matrix $B$ is symmetric of dimension $3$ with sums of columns $1.$ If $K_A$ has complement to \emph{$HW$ - matrix} then
also a matrix $K_B$ has this possibility. But it is impossible by Lemma \ref{oddlemma}.
In case (b) matrices $B,E$ are symmetric with sums of columns $3.$ Hence
$B = E = \left[
\begin{smallmatrix}
1 & 2\\
2 & 1
\end{smallmatrix}\right]$ and
\begin{equation}
p^{(2)}(A) = \left[
\begin{smallmatrix}
1 & 1 & 1 & 1\\
1 & 1 & 1 & 1\\
1 & 1 & 1 & 1\\
1 & 1 & 1 & 1
\end{smallmatrix}\right].
\end{equation}
In the matrix $p^{(2)}(K_A)$ we have rows of $1$ and $0$. They are linearly dependent from the rows of $p^{(2)}(A).$ Hence
$$\text{rk}p^{(2)}(K_A) = \text{rk}p^{(2)}(A) = 1$$
and
$$\text{d}(K_A) = 4-1 = 3 > 1.$$
From Lemma \ref{defect} the matrix $K_A$ has not complement to the \emph{$HW$-matrix}.
\vskip 2mm
\noindent
3. By the above points 1. and 2. we have that  $p^{(2)}(A)$ and  $p^{(3)}(A)$ are connected matrices.
As in the proof of Lemma \ref{oddlemma} using relations $\sim_{2}, \sim{_3}$ and conjugations of matrices we can assume that
$$[C\hskip 2mm A\hskip 2mm D] = [2\hskip 2mm A\hskip 2mm 3].$$
By assumption a sum of all rows of the above matrix has $1$ on a some position. We can see easily that it is impossible at the first and the third block.
For a matrix $A$ it is also impossible since $m$ is even.
This contradicts our assumption that $m < n.$
\vskip 2mm
\noindent
Let us assume that $m > 4.$ We shall consider three steps.
\vskip 1mm
\noindent
1. Assume that $p^{(2)}(A)$ is not connected. We have to consider two cases:
\vskip 2mm
\noindent
(a) $p^{(2)}(A)$ is a direct sum of two odd blocks,
\vskip 2mm
\noindent
(b) $p^{(2)}(A)$ is a direct sum of two even blocks.
\vskip 2mm
\noindent
Hence
$A = \left[
\begin{smallmatrix}
B & 3\\
3 & E
\end{smallmatrix}\right].$
In the case (a) since dimensions of $B,E$ are odd and sums of column of $A$ are $3$ we obtain that sums of column
of $B$ and $E$ are $0.$ Moreover, if $B$ is an odd diagonal submatrix of \emph{$HW$-matrix} then
by definition \ref{hwm} a sum of rows of $B$ should enclose $1.$ But this is impossible and also case (a) is impossible.
\vskip 2mm
\noindent
In case (b) since dimensions of $B,E$ are even and sums of column of $A$ are $3$ we obtain that sums of column
of $B$ and $E$ are $3.$ Moreover either the matrix $B$ or the matrix $E$ has rank $> 2.$  Assume the matrix $B$
has such a property.
If a matrix $K_A$ has complement, then a matrix $K_B$ has complement to \emph{$HW$-matrix}. But by induction it is impossible.
\vskip 2mm
\noindent
2. Assume that $p^{(3)}(A)$ is not connected. We have to consider two cases.
The same as in the step 1. 
\vskip 2mm
\noindent
(a) $p^{(3)}(A)$ is a direct sum of two odd blocks,
\vskip 2mm
\noindent
(b) $p^{(3)}(A)$ is a direct sum of two even bloks.
\vskip 1mm
\noindent
Hence
$A = \left[
\begin{smallmatrix}
B & 2\\
2 & E
\end{smallmatrix}\right].$
In the first case since dimensions of $B,E$ are odd and sums of column of $A$ are $3$ we obtain that sums of column
of $B$ and $E$ are $1.$ Moreover, either the matrix $B$ or the matrix $E$ has rank $> 2.$  Assume the matrix $B$ has such a property.
If a matrix $K_A$ has complement then (after permutation of indexes) a matrix $K_B$ has complement to \emph{$HW$-matrix}. 
But by Lemma \ref{oddlemma} it is impossible.
In the second case, since dimensions of $B,E$ are even and sums of column of $A$ are $3$ we obtain that sums of column
of $B$ and $E$ are $3.$ Moreover, either the matrix $B$  or $E$ has rank $> 2.$  Assume the matrix $B$ has such a property
If a matrix $K_A$ has complement then a matrix $K_B$ has complement to \emph{$HW$-matrix}. But by induction it is impossible
\vskip 2mm
\noindent
We can assume that matrices $p^{(2)}(A)$ and $p^{(3)}(A)$ are connected.
As in the previous cases we can assume that
$$[C\hskip 2mm A\hskip 2mm D] = [2\hskip 2mm A\hskip 2mm 3].$$
By definition \ref{hwm} a sum of all rows should enclose $1.$ Since $m$ is even and $m < n$ we have
a contradiction.
\vskip 2mm
\hskip 120mm $\Box$
\vskip 3mm
\section{Proof of the {\bf Main Lemma}}
We keep the notation from previous sections, but we also need a new definitions.
Denote by ${\cal P}_n$ an algebra of all subsets of the set $\{1,2,\dots,n\}.$ Let $\mid U\mid$ denote
the number of elements of a set $U\in {\cal P}_n$ modulo two.
We have an isomorphism of algebras $I:\F_{2}^{n}\to {\cal P}_n,$ 
where 
\begin{equation}\label{indicator}
I(x) = \{i\mid x_i = 1\},  x\in\F_{2}^{n} 
\end{equation}
is an indicator.
\vskip 1mm
\noindent

\begin{defin}\label{defj}
Let $A$ be a \emph{$HW$-matrix}. The function $J:{\cal P}_n\to {\cal P}_n$
is defined by
\begin{equation}\label{formulaJ}
J(U) = \{s\mid \sum_{i\in U} A_{is} = 1\},
\end{equation}
where $U\in {\cal P}_n.$
\end{defin}
\begin{rem}\label{examplej}
In what follows we shall use a formula {\em (\ref{bformula})} with
a basis $b_{i}, 1\leq i\leq n-1.$
Let us consider a map $l:{\cal P}_{n}\to\F_2[x_1,\dots,x_{n-1}]$
given by a formula
\begin{equation}\label{chara}
l_{Z}: = \sum_{i\in Z}x_{i}.
\end{equation}
In this language the formula {\em (\ref{bformula})} for $k = 2,3$ we can write as
$$\sum_{j=1}^{n-1}p^{(k)}A_{ji}x_{j} = l_{S}$$
where $S = \{p^{(k)}(A_{1,i}),p^{(k)}(A_{2,i}),...,p^{(k)}(A_{n-1,i})\}.$
\end{rem}

\vskip 2mm
\begin{prop}\label{propJ}
The map $J$ has the following properties:
\begin{enumerate}
\item
$U\neq 0,1$ then $J(U)\neq 0,$ here $0,1$ denote the trivial additive and multiplicative element of the algebra ${\cal P}_n$ respectively;
\item
$J(U+1) = J(U)$ where $U+1 = U'$ denotes a complement of the subset $U$ in the set $\{1,2,...,n\};$
\item
$J(\{i\}) = \{i\}, i=1,2,...,n;$
\item
if $\mid U\mid = 1$ then $J(U)\subset U;$
\item
if $\mid U\mid = 0$ then $J(U)\subset U'.$
\end{enumerate}
\end{prop}
{\bf Proof:} Elementary calculations with support of the matrix (\ref{matrix1}).
\vskip 2mm
\hskip 120mm $\Box$
\vskip 3mm
\noindent
Any polynomial of $\F_2[x_1,x_2,\dots,x_n]$ we shall identify with a polynomial map $\F_{2}^{n}\to\F_2.$
Hence by indicator function (\ref{indicator}) the formula (\ref{chara}) has the following presentation $l_{Z}(e_j) = \{j\in Z\},$
where $Z\in {\cal P}_n.$
Since the transgressive elements $T_i\in\F_2[x_1,\dots,x_{n-1}]$ we define
a split monomorphism of rings $\F_2[x_1,\dots,x_{n-1}]\stackrel{\phi}{\to}\F_2[x_1,\dots,x_n]$ 
such that $\bar{T}_i = \phi(T_i)\in \F_2[x_1,\dots,x_n], i = 1,\dots,n.$
Here, $\phi(x_i) = x_i + x_n, i = 1,2,\dots,n-1.$
Obviously $\# D = \#\phi(D).$
 
\vskip 1mm
\noindent
From definition, for polynomial functions $\bar{T}_i$ we have  
$\bar{T}_{i}(e_j) = \delta_{ij},$
where $1\leq i,j\leq n$ and $e_i\in(\F_2)^{n}$ is the standard basis.
Hence, by the isomorphism (\ref{indicator}) a map $J$ (see Definition \ref{defj}) is equivalent to a function
$T:\F_{2}^{n}\to\F_{2}^{n}, T(x) = (\bar{T}_1(x),\bar{T}_2(x),\dots,\bar{T}_n(x)),$
where $x\in\F_{2}^{n}.$
Hence and from an equation (\ref{elemmaa}) we have a commutative diagram
\begin{equation}\label{diag1}
\begin{diagram}
\node{\F_2^n}\arrow{s,r}{I}\arrow{e,t}{T}\node{\F_2^n}\arrow{s,r}{I}\\
\node{{\cal P}_n}\arrow{e,t}{J}\node{{\cal P}_n}
\end{diagram}.
\end{equation}

We shall use these observations
in the proof of the {\bf Main Lemma}. 
Moreover, we shall apply a remark that 
homogeneous polynomials of degree $2$ are recognized by their polynomial functions.
Let $S, Z_1, Z_2\in {\cal P}_{n}.$ From definition if 
$$\sum_{i\in S}\bar{T_i} = l_{Z_1}\cdot l_{Z_2}$$
then $S = Z_1\cap Z_2.$ 
\vskip 3mm
\begin{prop}\label{prop31}
The following conditions are equivalent.
\vskip 1mm
{\em (i)} $\sum_{i\in S}\bar{T}_i = l_{Z_1}\cdot l_{Z_2}$
\vskip 1mm
{\em (ii)} $\forall U\in {\cal P}_{n} |J(U)S| = |U Z_1|\cdot |U Z_2|.$
\end{prop}
\vskip 3mm
\noindent
\noindent
{\bf Proof:} 
We shall use (\ref{diag1}) and an isomorphism $I.$
Let $x\in\F_2^n, U = I(x).$ 
We have 
$$\sum_{i\in S}\bar{T}_i(x) = \sum_{i\in S\cap I(\bar{T}(x))} 1 = \mid I(\bar{T}(x))\cap S\mid = \mid J(I(x))\cap S\mid = \mid J(U)\cap S\mid.$$ 
From the other side $$l_{Z_1}(x)\cdot l_{Z_2}(x) = \sum_{i\in Z_1\cap I(x))} 1\cdot \sum_{i\in Z_2\cap I(x))} 1 = |U Z_1|\cdot |U Z_2|.$$
This finishes a proof.
\vskip 2mm
\hskip 120mm $\Box$
\vskip 3mm
\newtheorem{corol}{Corollary}
\begin{corol}\label{corprop3}
Let us assume the condition {\em (ii)} of Proposition \ref{prop31}, then
\vskip 1mm
\noindent
$1$. $|Z_1|$ or $|Z_2|$ is even,
\vskip 1mm
\noindent
$2$. if $S\neq 0$ then $|Z_1|$ and $|Z_2|$ are even
\vskip 1mm
\noindent
$3$. if $S\neq Z_1$ and $S\neq Z_2$ then $Z_1\cup Z_2 = 1.$
\end{corol}
{\bf Proof:} 
$1$. Since $J(1) = J(\{1,2,\dots,n\}) = 0$ the condition is true.
\vskip 1mm
\noindent
$2$. Since $J(U) = J(U') = J(U+1)$ we have
$$|UZ_{1}||UZ_{2}| = |(1+U)Z_{1}||(1+U)Z_{2}|.$$
Hence
$$|Z_{1}||Z_{2}| + |Z_{1}||UZ_{2}| + |Z_{2}||UZ_{1}| = 0.$$
From a point 1. we can assume that $|Z_1| = 0$ (or $|Z_2| = 0$) and $|Z_{2}||UZ_{1}| = 0.$
If $|Z_{2}| = 1$ then $\forall U\in {\cal P}_n, |UZ_1| = 0$ and $Z_1 = 0.$ 
Since $S = Z_1\cdot Z_2\neq 0$ we have a contradiction.
\vskip 1mm
\noindent
$3$. Let $a\in Z_1\setminus S, b\in Z_2\setminus S$ and $c\notin Z_1\cup Z_2.$
Put $U = \{a,b,c\}.$ We have $J(U)S\subset US = 0$ and $UZ_1 = \{a\}, UZ_2 = \{b\}.$
Hence
$$0 = |J(U)S| = |UZ_1||UZ_2| = 1\cdot 1 = 1.$$
This is a contradiction.
\vskip 2mm
\hskip 120mm $\Box$
\vskip 3mm
\begin{defin}
Define
$$\sigma^{S}_{a} := \sum_{i\in S} A_{a,i},$$
where  $a\in \{1,2,\dots,n\}, S\subset \{1,2,\dots,n\}$ and $A\in {\cal D}^{n\times n}.$
\end{defin}
Let us present relations between the above definition and the function $J.$ 
\begin{prop}\label{rowsum}
Let $A$ be $(n\times n)$ \emph{$HW$-matrix}, $a, b\in\{1,2,\dots,n\}$ and $S\in {\cal P}_n.$
Then
\vskip 1mm
\noindent
1. $|J(\{a,b\})S| = \sigma^{S}_{a} + \sigma^{S}_{b},$ where $a,b \notin S$;
\vskip 1mm
\noindent
2. $|J(\{a,b\})S| = \sigma^{S}_{a} + \sigma^{S}_{b} + A_{a,b} + 1,$ where $a\notin S , b\in S$; 
\vskip 1mm
\noindent
3.  $|J(\{a,b\})S| = \sigma^{S}_{a} + \sigma^{S}_{b} + A_{a,b} + A_{b,a},$ where $a, b\in S.$
\end{prop}
{\bf Proof:} 1. By a point 5. of Proposition \ref{propJ} we know that $J(\{a,b\})\subset \{a,b\}'.$ If $J(\{a,b\})S = \emptyset$  we are done.
On the contrary we shall consider the following cases. 
\vskip 1mm
\noindent
(a) Assume $|S| = 1$ and $|J(\{a,b\})S| = 1.$ We have two rows, which correspond to $a$ and $b,$
\begin{equation}\label{2rows}
\begin{array}{ccccccccc}
2&2&\dots& 2& 2\\
2&3&\dots&3&2
\end{array}
\end{equation}
with a number of columns equal to $|S|,$ and a number of columns with different coefficients equal to $J(\{a,b\}).$ 
Hence a sum of the upper row is equal to $2$ and a sum of the down row is equal to $3.$ This finishes a proof in this case.
\vskip 1mm
\noindent
(a') Assume $|S| = 1$ and $|J(\{a,b\})S| = 0.$ We also have (\ref{2rows}) and a sum of the upper row is equal to $2$ and a sum of the down row is also equal to $2.$ 
This finishes a proof in this case.
\vskip 1mm
\noindent
(b) Assume $|S| = 0.$ Then again we have two subcases $|J(\{a,b\})S| = 1,$ then a sum of the upper row of (\ref{2rows}) is equal to $0$ and a sum of the down row is equal to $1.$
The proof of the case is complete.
When $|J(\{a,b\})S| = 0$ a sum of the upper row of (\ref{2rows}) is $0$ and a sum of the down row is also $0.$
This finished a proof of point 1. The proofs of other cases are similar and we put it as an exercise.
\vskip 2mm
\hskip 120mm $\Box$
\vskip 3mm
Using the above language we shall prove that for a \emph{$HW$-manifold} there exists only a limited number of trangressive elements which are a product of degree one nontrivial polynomials.
\begin{prop}\label{trans1}
Let $A$ be a $(n\times n)$ \emph{$HW$-matrix}, ($n > 3$) then there does not exist not empty set $S\subset\{1,2,\dots,n\}$ such that
\begin{equation}\label{rozk1}
\forall U\in {\cal P}_n |J(U)S| = |US|.
\end{equation}
\end{prop}
{\bf Proof:}
It is the case $S = Z_1 = Z_2.$ Let us assume (\ref{rozk1}). 
We are going to divide the proof into four steps.
\vskip 1mm
\noindent
{\em Step 1.}
We claim that, if $a_1, a_2\notin S$ and $b\in S$ then  $A_{a_{1},b} = A_{a_{2},b}.$
In fact, from (\ref{rozk1}) for $U = \{a_1, a_2\},$ $|J(\{a_1, a_2\})S| = |\{a_1, a_2\}S| = 0.$  
By Proposition \ref{rowsum} (1.), $\sigma^{S}_{a_{1}} = \sigma^{S}_{a_{2}} := \sigma.$
If $a\notin S$ then from Proposition \ref{rowsum} (2.) 
$$1 = |J(\{a,b\})S| = |\{a,b\}S| = \sigma_{a}^{S} + \sigma_{b}^{S} + A_{a,b} + 1 =  \sigma + \sigma_{b}^{S} + A_{a,b} + 1.$$
Hence $\forall a\notin S, A_{a,b} = \sigma + \sigma_{b}^{S}$ and {\em Step 1.} is proved.
\vskip 1mm
\noindent
{\em Step 2.}
We claim that, if $US = 0$ then $J(U)S = 0.$  In fact from {\em Step 1.} all elements (numbers of columns) of $J(U)$ which are considered have not the first indexes from $S$ and
are equal each other. Then $J(U)S = 0.$ 
\vskip 1mm
\noindent
{\em Step 3.}
We claim that, if $S\neq 0$ then $\# S = n-1.$
From {\em Step 2.} if $0\neq U\subset S'$ then $J(U)S'\neq 0.$
Let $B$ be a diagonal submatrix of the matrix $A$ related to the set $S'.$
Then $B$ is a quadratic matrix with $1$ on the diagonal and $2,3$ otherwise.
Moreover all sums of rows of $B$ have at some position an element $1.$ Hence, the only possible matrix $B$ is $(1\times 1)$ matrix.
\vskip 1mm
\noindent 
{\em Step 4.} We claim that, if $S\neq 0$ then $n\leq 3.$
For the proof, let us assume that $n > 3.$ From the {\em Step 3.} we can assume that $S = \{2,3,\dots,n\}.$
Let $l_2$ denote a number of $2$ at the first column of $A.$ We shall prove that $|l_2| = 0.$
In fact, we can assume that $0 < l_2 < n-1$ and at the first column, from the top we have first $2$ then going down we have $3.$
On the contrary, suppose that $l_2$ is odd and let $v$ be a sum of the first $2l_2 +1$ rows.
Since $l_2 + 1$ is even $v$ has not $1$ on places $1,2,\dots,l_2 + 1.$ Then it has $1$
on the position $> l_2 + 1.$  Hence there exists $k\geq l_2 + 1$ such that $A_{1,r}\neq A_{k,r}$
or equivalently $A_{1,r} +  A_{k,r} = 1.$ Let us consider a diagonal submatrix 
\begin{equation}
\left[
\begin{smallmatrix}
1 & \ast & A_{1,r}\\
2 & 1 & A_{k,r}\\
3 & \ast & 1
\end{smallmatrix}\right].
\end{equation}
A sum of elements at the first column and at the third column is $0$, then it at the second column has to be $\neq 0.$
Let $U = \{1,k,r\}.$ Since $j(U)\subset U$ and $n > 3$, $J(U) = \{k\}.$ Finally
$$1 = |\{k\}S| = |J(U)S| = |US| = |\{k,r\}S| = 0.$$
That is a contradiction and $l_2$ is even. Moreover if $l_3$ is a number of $3$ at the first column
then $\mid l_3 = n-1 - l_2\mid = 0$ and a sum $1 + l_2\ast 2 + l_3\ast 2 = 1$. But a sum of all rows
is zero and we have a contradiction. This finishes a proof.
\vskip 2mm
\hskip 120mm $\Box$
\vskip 3mm
\begin{corol}\label{squar}
At the space Im($d_2$) we have not squares.
\end{corol}
{\bf Proof:} If $l_{Z}\in$ Im($d_2$), then $S = Z = Z.$ For $n > 3$ it is impossible.
\vskip 2mm
\hskip 120mm $\Box$
\vskip 3mm
\begin{prop}\label{trans2}
Let $S,Z\subset\{1,2,\dots,n\}$ such that $0\neq S\neq Z.$
Let $A, J$ be as in Proposition \ref{trans1}. 
Assume that
$$\forall U\in {\cal P}_n \mid J(U)S\mid = \mid US\mid\cdot\mid UZ\mid$$
then $\# S = 2, \mid Z\mid = 0$ and $S\subset Z.$ 
\end{prop}
{\bf Proof:} On the beginning we claim that up to permutation and conjugation,
\begin{equation}
A = \left[
\begin{smallmatrix}
\ast & 2 & \ast\\
\ast & B & \ast\\
\ast & 3 & \ast
\end{smallmatrix}\right],
\end{equation}
where $B$ is a symmetric matrix with a column sums $3.$ Moreover a block $2$ has rows indexed by numbers 
from the set $Z\setminus S$
and a block $3$ has rows indexed by numbers from the set $1+Z = Z'.$
In fact, from Proposition \ref{prop31}, $S\subset Z$ and Corollary \ref{corprop3}, $S\subset Z$ and $|S| = |Z| = 0.$
Let us change the indexes of $A$ such that  
\begin{equation}
A = \left[
\begin{smallmatrix}
\ast & E & \ast\\
\ast & B & \ast\\
\ast & F & \ast
\end{smallmatrix}\right],
\end{equation}
and $E$ has rows indexed by numbers from the set  $Z\setminus S,$ $B$ 
has rows indexed by numbers from $S$ and $F$ is indexed by $1+Z = Z'.$
From the point 1 of Proposition \ref{rowsum}, for $a,b\notin S$
$$\sigma_{a}^{S} + \sigma_{b}^{S} = |J(\{a,b\}S| = |\{a,b\}S|\cdot |\{a,b\}Z| = 0.$$
Hence $\sigma_{a}^{S} = \sigma_{b}^{S}.$ Let $\sigma := \sigma_{a}^{S},$ for $a\notin S.$
\vskip 1mm
\noindent
By the point 2 of Proposition \ref{rowsum} for $b\in S$ and $a\notin Z,$
\begin{equation}\label{point2}
A_{a,b} = \sigma + \sigma_{b}^{S}.
\end{equation}
From the above all columns of the matrix $F$ are constant.
Again from the point 2 of Proposition \ref{rowsum} for $b\in S, a\in Z\setminus S,$
\begin{equation}\label{point4}
A_{a,b} = \sigma + \sigma_{b}^{S} + 1.
\end{equation}
It follows that also columns of the matrix $E$ are constant.
Let us conjugate columns of the matrix $A$ such that $E = 2.$
In that case $\sigma = 0$ since for $a\in Z\setminus S$ we have $\sigma = \sigma_{a}^{S} = |S|\cdot 2 = 0.$
\vskip 2mm
\noindent
From (\ref{point4}), for $b\in S, 2 = 0 + \sigma_{b}^{S} + 1.$ Hence $\sigma_{b}^{S} = 3$
and $F = 3,$ because from the formula (\ref{point2}) $A_{a,b} = 0 + 3,$ for $a\in Z'\hskip 2mm \text{and}\hskip 2mm b\in S.$
Finally, from Proposition \ref{rowsum} for $a,b\in S$ we have
$$A_{a,b} + A_{b,a} = 3 + 3 + A_{a,b} + A_{b,a} = \sigma_{a}^{S} + \sigma_{b}^{S} + A_{a,b} + A_{b,a} =$$
\begin{equation}
= |J(\{a,b\})S| = |\{a,b\}S|\cdot |\{a,b\}Z| = 0.  
\end{equation}
To finish a proof it suffices to apply Lemma \ref{evenlemma}.
\vskip 2mm
\hskip 120mm $\Box$
\vskip 3mm
\begin{prop}\label{trans3}
We keep the notation from the previous propositions.
Let us assume $S, Z_1, Z_2\in {\cal P}_n$ such that $0\neq S, S\neq Z_1, S\neq Z_2$ and
$$\forall U\in {\cal P}_n |J(U)S| = |UZ_1|\cdot |UZ_2|$$
then $\# S = 1, |Z_1| = |Z_2| = 0$ and $Z_1 + Z_2 = 1.$
\end{prop}
{\bf Proof:} A proof is similar to the proof of Proposition \ref{trans2}.
On the beginning we show that (up to permutation and conjugation)
\begin{equation}
A = \left[
\begin{smallmatrix}
\ast & 2 & \ast\\
\ast & B & \ast\\
\ast & 3 & \ast
\end{smallmatrix}\right],
\end{equation}
where $B$ is a symmetric matrix of odd dimension with sums of columns $1,$ a block $2$
is indexed by the set $Z_1\setminus S$ and a block $3$ is indexed by the set $Z_2\setminus S.$
In fact, from assumption and Corolarry \ref{corprop3}, $S = Z_{1}Z_{2}, |Z_{1}| = |Z_{2}| = 0$ and
$Z_1 + Z_2 = 1.$ Hence $|S| = 1.$ Let us change the order of rows in the matrix $A$ such that
\begin{equation}
A = \left[
\begin{smallmatrix}
\ast & E & \ast\\
\ast & B & \ast\\
\ast & F & \ast
\end{smallmatrix}\right]
\end{equation}
and $E$ is indexed by $Z_1\setminus S,$ $B$ by $S$ and $F$ by $Z_2\setminus S.$
From Proposition \ref{rowsum} we have $\forall a,b\in Z_1\setminus S, \sigma_{a}^{S} = \sigma_{b}^{S} :=\sigma_{E}.$
With similar consideration we have $\forall a,b\in Z_2\setminus S, \sigma_{a}^{S} = \sigma_{b}^{S} := \sigma_{F}.$
Moreover, by Proposition \ref{rowsum} (2) for $b\in S$ and $a\in Z_1\setminus S,$
\begin{equation}\label{pointthree}
A_{a,b} = \sigma_{E} + \sigma_{b}^{S} +1.
\end{equation}
From the above, all columns of the matrix $E$ are the same.
By analogy for $b\in S$ and $a\in Z_2\setminus S,$
\begin{equation}\label{pointfive} 
A_{a,b} = \sigma_{F} + \sigma_{b}^{S} + 1
\end{equation}
and columns of the matrix $F$ are also constant.
Let us conjugate columns of $A$ such that $E = 2.$
Then $\sigma_{E} = 2,$ because for $a\in Z_1\setminus S,\sigma_E = \sigma_{a}^{S} = |S|\cdot 2 = 2$
and for $b\in S, \sigma_{b}^{S} = 1.$ The last equality follows from (\ref{pointthree}) because
$2 = 2 + \sigma_{b}^{S} + 1.$
Similarly, by (\ref{pointfive}) for $b\in S$ and $a\in Z_2\setminus S,$
we have $A_{a,b} = \sigma_{F} + 1 + 1 =\sigma_F$ and the matrix $F$ is constant and equal to $\sigma_{F}.$
Finally, a matrix $B$ is symmetric since from Proposition \ref{rowsum}
$$\sigma_{a}^{S} + \sigma_{b}^{S} + A_{a,b} + A_{b,a} = |J(\{a,b\}S|$$
what means,
$$ A_{a,b} + A_{b,a} = |\{a,b\}Z_{1}|\cdot |\{a,b\}Z_{2}| = 0.$$
We have still to show that $\sigma_F = 3.$
In fact from assumption a column's sums of $B$ are $1.$ 
Since $B$ is symmetric this same is true for rows.
Let us calculate a sum of some column of $A:$
$$(|Z_{1}| - |S|)2 + 1 + (|Z_{2}| - |S|)\sigma_{F} = 2 + 1 + \sigma_F = 3 + \sigma_F = 0.$$
To finish a proof of Proposition we have to apply Lemma \ref{oddlemma}.
\vskip 2mm
\hskip 120mm $\Box$
\vskip 3mm
Summing up we have the following two possibilities:
\vskip 2mm
\noindent
I. $\# S = 1$ and $Z_1 + Z_2 = 1;$
\vskip 1mm
\noindent
II. $\# S = 2$ and $S = Z_1, S\neq Z_2$ or $S = Z_2. S\neq Z_1.$
\vskip 3mm
\noindent
Let us recall that Im($d_2$) is a $n$-dimensional $\Z_2$- space generated by $T_i, i = 1,2,3,\dots,n.$                
We are interested in description of the set $D$ of 
elements in Im$d_2$ which are a product of two nontrivial linear polynomials, see (\ref{Span}).
We claim that $D\leq n+1.$
In what follows, if it does not give a contradiction we shall write $T_i$ for $\bar{T}_i, i = 1,\dots,n.$
\begin{lem}
Let $w\in D,$ then $w = T_i$ or $w = T_j + T_k$ for some $1\leq i,j,k\leq n.$
\end{lem}
{\bf Proof:} 
On the beginning we shall prove that $T_i + T_j$ is a product of two nontrivial linear polynomials if and only if
$T_i, T_j$ have a common component. It means there exists ${\rm p} \neq 0$ s.t. ${\rm p}|T_i$ and ${\rm p}|T_j.$
Let $T_i + T_j$ have a common component, then from the above case II we can assume that 
$j = i + 1$ and the matrix $A$ enclose:
$$...$$
$$2\hskip 3mm 2$$
$$
\left[
\begin{array}[c]{ccc}
1 & 2\\
2 & 1
\end{array}\right]$$
$$3\hskip 3mm 3$$
$$...\hskip 1mm.$$
By definition
$$T_i = (x_1 +\dots + x_i + x_{i+1})(x_i + x_{i+2} +\dots + x_n)$$
$$T_{i+1} = (x_1 +\dots + x_i + x_{i+1})(x_{i+1} +\dots + x_n).$$
For the proof of the opposite conclusion we shall need
\begin{defin}
Let $X$ be a subset of some monoid. By $\Gamma_X$ we define a graph with the vertex
set $X$ and two vertices $a,b$ are connected by an edge $a\stackrel{f}{\rule{25pt}{0.4pt}} b$ if and only if $f | a$ and $f | b.$
Put $\Gamma := \Gamma_{T_1,T_2,\dots, T_n}.$
\end{defin}
We claim that for $n > 3$ the graph
\begin{equation}\label{graphcom}
i\stackrel{f}{\rule{25pt}{0.4pt}} j\stackrel{g}{\rule{25pt}{0.4pt}}k
\end{equation}
is not a subgraph of $\Gamma,$ where $i:= T_i, i = 1,2,\dots,n.$
In fact we have two possibilities:
\vskip 2mm
\noindent
1. $f = g.$ Let $i = 1, j = 2, k = 3$ and let $\mathfrak{I}$ be an ideal generated by $(f,T_4,\dots,T_n)$ in the polynomial ring. 
Since there exist a nontrivial solution of system of $(n-2)$ linear equation in $(n-1)$ linear space an algebraic set V($\mathfrak{I}$)  is not trivial.
It means $0\neq x\in$ V($\mathfrak{I}$). From definition $x\in$ V($\mathfrak{I}'$), where $\mathfrak{I'}$ is an ideal generated by $(T_1,T_2,\dots, T_n).$
But it is impossible.
\vskip 2mm
\noindent
2. $f\neq g.$  Using permutation of indexes and conjugation we can assume that in \emph{$HW$-matrix} $A, j = i + 1, k = i + 2.$
Recall that $S = \{i, i + 1\}$ and $A$ is as in Lemma \ref{evenlemma}. Hence it has a diagonal block related to rows (columns) $\{i, i + 1, i + 2\}$ 
\begin{equation}\label{star}
\left[
\begin{smallmatrix}
1 & 2 & b\\
2 & 1 & a\\
3 & 3 & 1
\end{smallmatrix}\right],
\end{equation}
and a matrix $A$ has upper two first columns of $(\ref{star})$ only elements $2$, but lower only elements $3.$
Let us consider polynomials $T_i, T_{i+1}$ and $T_{i+2}$ for $x_s = 0, s\notin \{i, i + 1, i + 2\}$ and denote it by $\hat{T}_i$
respectively.
We have
$$\hat{T}_{i} = (x_i + x_{i + 1})(x_i + x_{i + 2})$$
and
$$\hat{T}_{i + 1} = (x_i + x_{i + 1})(x_{i + 1} + x_{i + 2}).$$
The both polynomials are divided by $(x_i + x_{i + 1}).$ 
Hence  $\hat{T}_{i + 1}$ and  $\hat{T}_{i + 2}$ are divided by $(x_{i + 1} + x_{i + 2}).$
From the above we can observe that 
\begin{equation}\label{fgt}
\hat{T}_{i + 2} = (x_{i + 1} + x_{i + 2})(x_{i + 2} + x_{i}).
\end{equation}
By (\ref{fgt}) and definition we get $a\neq b.$ Hence
a sum of all columns of the matrix ($\ref{star}$) are equal to $0.$ But it is impossible, since $n > 3.$
This finishes a proof of our claim and we have
\vskip 1mm
\begin{corol}
For $n> 3$ all connected components of a graph $\Gamma$ are points or edges
$i\stackrel{f}{\rule{25pt}{0.4pt}} j.$
\end{corol}
\vskip 2mm
\hskip 120mm $\Box$
\vskip 3mm
\begin{corol}
For $n > 3, D = \{T_1, T_2,..., T_n\}$ or {\em D} =  $\{T_1, T_2,..., T_n, T_i + T_j\}$ for some $1\leq i,j\leq n.$
\end{corol}
{\bf Proof:}
Conversely, suppose that edges 
$$1\stackrel{f}{\rule{25pt}{0.4pt}} 2 \hskip 2mm\text{and}\hskip 2mm  3\stackrel{g}{\rule{25pt}{0.4pt}} 4$$ are components of the graph $\Gamma.$
Let us consider an ideal  $\mathfrak{J}$ = $(f,g, T_5,\dots, T_n)$
in polynomial ring.
Since there exist a nontrivial solution of system of $(n-2)$ linear equation in $(n-1)$ linear space an algebraic set V($\mathfrak{J}$)  is not trivial.
It means $0\neq x\in$ V($\mathfrak{J}$). But from definition $x\in$ V($\mathfrak{I}$'), where  $\mathfrak{J}'$ is an ideal generated by $(T_1,T_2,\dots, T_n).$
But it is impossible. This finishes a proof.
\vskip 2mm
\hskip 120mm $\Box$
\vskip 3mm
Let us prove the {\bf Main Lemma}.
\vskip 5mm
\noindent
{\bf Main Lemma}
{\em Let $n > 3$, then there are the following possibilities for the structure of the set $D$:
\vskip 1mm
1. $D = \{T_1,T_2,\dots, T_n\};$
\vskip 1mm
2. $D = \{T_1,T_2,\dots, T_n, T_i + T_j\},$ and we can find a polynomial ${\rm p}$ of degree one such that
${\rm p}\mid T_i$ and ${\rm p}\mid T_j$ for some $1\leq i,j\leq n.$
In the second case we can rediscover the set of generators $T_1, T_2,\dots, T_n.$}
\vskip 3mm
\noindent
{\bf Proof:}
We start from the simple observation.
If $i\neq j$ and $T_i = {\rm p\cdot q}, T_j = {\rm p\cdot r}$ then $\forall i = 1,2,\dots,n\hskip 3mm {\rm q + r}$ is not divided $T_i.$ 
In fact ${\rm q + r}\neq {\rm p}$ since in other case $T_i + T_j = {\rm p(q + r) = p^2}.$ 
By Corollary \ref{squar} it is impossible.
Hence $T_i$ and $T_j$ are also not divided by ${\rm q + r}.$ Moreover, if $T_r = {\rm (q + r)s}$
then $T_i + T_j + T_r = {\rm (q + r)(p + s)}.$ By Proposition \ref{trans2}, a decomposition for $\# S = 3$ is impossible.
Let us prove the second point of the above lemma. From definition the graph
$\Gamma_{T_1,...,T_n}$ has connected components which are vertices for $r\notin\{i,j\},$
of the triangle with vertices $T_i, T_j, T_i + T_j$ and a constant label which is a component
of $T_i$ and $T_j.$
Let $T_i = {\rm p\cdot q}$ and $T_j = {\rm p\cdot r}$ then $T_i + T_j = {\rm p(q + r)}.$
The triangle is a connected component of a graph because by (\ref{graphcom}) for $r\notin\{i,j\}$ 
elements ${\rm p,q,r}$ do not divide $T_r.$ Also from the above simple observation, the element ${\rm (q + r)}$ does not divide $T_r.$
\vskip 2mm
We continue the proof of the {\bf Main Lemma}. Let $w = \xi\eta$ where $\xi$ and $\eta$ are linear polynomials.
Let us define $s(w) := \xi + \eta.$
Since \emph{$HW$-manifolds} are oriented $\sum_{i}s(T_i) = 0.$
We claim that if
$T_i + T_j\in D,$ then $s(\xi) + s(\eta)$ recognizes subsets of order two of the set $\{T_i, T_j, T_i + T_j\}.$
In fact, let $T_i = {\rm p\cdot q}, T_j = {\rm p\cdot r},$ then $T_i + T_j = {\rm p(q + r)}$ and
$s(T_i) + s(T_j) = {\rm q+r}, s(T_i) + s(T_i + T_j) = {\rm r}, s(T_j) + s(T_i + T_j) = {\rm q}.$
\vskip 2mm
Let $n > 3$, then there are the following possibilities for the structure of the set $D$:
\vskip 1mm
1. $D = \{T_1,T_2,\dots, T_n\};$
\vskip 1mm
2. $D = \{T_1,T_2,\dots, T_n, T_i + T_j\},$ for some $1\leq i,j\leq n.$
Let $n > 3$ if $D$ has $n$ elements we are done.
If it has $(n + 1)$ elements then the graph $\Gamma_{T_1,T_2,...,T_n}$ 
has $(n-2)$ discrete connected components $D^{c}$ and a triangle.
We proceed in two steps:
\vskip 1mm
\noindent
1. Put $s_{D^{c}} := \sum_{a\in D^{c}}s(a)$
\vskip 1mm
\noindent
2. From the triangle we take a unique pair $\xi,\eta$ such that
$$s(\xi) + s(\eta) + s_{D^{c}} = 0.$$
Hence $\{T_1, T_2,\dots, T_n\} = \{\xi,\eta\}\cup D.$
This finishes a proof of the {\bf Main Lemma}.
\vskip 2mm
\hskip 120mm $\Box$
\vskip 3mm
For illustartion of possibilities of the structure of the set $D$ we present two examples.
\begin{examp}
Let $G\subset {\cal D}^5$ correspond to \emph{$HW$-matrix}   
$\left[
\begin{smallmatrix}
1 & 2 & 2 & 2 & 2\\
2 & 1 & 3 & 2 & 2\\
3 & 2 & 1 & 3 & 2\\
3 & 2 & 3 & 1 & 3\\
3 & 3 & 3 & 2 & 1
\end{smallmatrix}\right].$
 
The set 
\begin{equation}
\begin{split}
D =& \{T_1 = (x_1+x_2)(x_1+x_3+x_4), T_2 = (x_1+x_2+x_3+x_4)x_2,\\ 
& T_3 = (x_1+x_3)(x_2+x_3+x_4)), T_4 = (x_1+x_2+x_4)(x_3+x_4),\\
& T_5 = (x_1+x_2+x_3)x_4\}.
\end{split}
\end{equation}
\end{examp}
\vskip 2mm
From Remark \ref{remark1} the above group is isomorphic to the group $\Gamma_1$ of the
example \ref{example1}.
The next example illustrate the second case of the {\bf Main Lemma}.
\begin{examp} Let a matrix
$\left[
\begin{smallmatrix}
1 & 2 & 2 & 2 & 2\\
2 & 1 & 3 & 2 & 2\\
2 & 2 & 1 & 3 & 2\\
2 & 2 & 2 & 1 & 3\\
3 & 3 & 2 & 2 & 1
\end{smallmatrix}\right] \in {\cal D}^{5\times 5}$ be the second \emph{$HW$-matrix} of dimension $5.$

In this cases we have
\begin{equation}
\begin{split}
D = & \{T_1 = (x_1+x_2+x_3+x_4)x_1, T_2 = (x_1+x_2+x_3+x_4)x_2,\\ 
& T_3 = (x_1+x_3+x_4)(x_2+x_3), T_4 = (x_1+x_2+x_4)(x_3+x_4),\\
& T_5 = (x_1+x_2+x_3)x_4, T_1 + T_2\}.
\end{split}
\end{equation}
\end{examp}
{\bf Acknowledgments:} Authors would like to thank the referee for many constructive remarks,
which improve the exposition. The second author would like to thank R. Lutowski for help
during preparing the manuscript.

\end{document}